\documentclass[12pt]{article}
\usepackage{amsmath,amsthm,amssymb,hyperref,setspace,geometry}
\usepackage[all]{hypcap}
\usepackage{mathrsfs}
\title{On $c$-embedded subgroups of finite groups}
\author{Julian Kaspczyk\footnote{
		Institut für Algebra, Technische Universität Dresden, 01069 Dresden, Germany, \newline Email address: julian.kaspczyk@gmail.com}
}
\date{}
\theoremstyle{definition}

\newtheorem{innercustomthm}{Theorem}
\newenvironment{customthm}[1]
{\renewcommand\theinnercustomthm{#1}\innercustomthm}
{\endinnercustomthm}

\newtheorem{innercustomdef}{Definition}
\newenvironment{customdef}[1]
{\renewcommand\theinnercustomdef{#1}\innercustomdef}
{\endinnercustomdef}

\newtheorem{innercustomlem}{Lemma}
\newenvironment{customlem}[1]
{\renewcommand\theinnercustomlem{#1}\innercustomlem}
{\endinnercustomlem}

\newtheorem{innercustomex}{Example}
\newenvironment{customex}[1]
{\renewcommand\theinnercustomex{#1}\innercustomex}
{\endinnercustomex}

\newtheorem{innercustomqu}{Question}
\newenvironment{customqu}[1]
{\renewcommand\theinnercustomqu{#1}\innercustomqu}
{\endinnercustomqu}

\expandafter\let\expandafter\oldproof\csname\string\proof\endcsname
\let\oldendproof\endproof
\renewenvironment{proof}[1][\proofname]{%
	\oldproof[\bfseries\scshape #1]%
}{\oldendproof}

\usepackage{stackengine,scalerel}
\stackMath

\renewcommand{\le}{\leqslant}
\renewcommand{\ge}{\geqslant}

\begin{document}
	\maketitle
	\begin{abstract}
		Let $G$ be a group and $H \le K \le G$. We say that $H$ is \textit{$c$-embedded} in $G$ with respect to $K$ if there is a subgroup $B$ of $G$ such that $G = HB$ and $H \cap B \le Z(K)$. Given a finite group $G$, a prime number $p$ and a Sylow $p$-subgroup $P$ of $G$, we investigate the structure of $G$ under the assumption that $N_G(P)$ is $p$-supersolvable or $p$-nilpotent and that certain cyclic subgroups of $P$ with order $p$ or $4$ are $c$-embedded in $G$ with respect to $P$. New characterizations of $p$-supersolvability and $p$-nilpotence of finite groups will be obtained.
		
		\medskip
		
	\end{abstract}
	\section{Introduction}
	All groups in this paper are implicitly assumed to be finite. We use standard notation and terminology, see for example \cite{Huppert} or \cite{Isaacs}. Throughout, $p$ denotes an arbitrary but fixed prime.
	
	Recall that a group $G$ is said to be \textit{$p$-nilpotent} if $G$ has a normal Hall $p'$-subgroup. This concept plays an important role in local finite group theory, and many criteria for $p$-nilpotence of finite groups can be found in the literature. For example, a well-known theorem of Burnside asserts that a group $G$ with Sylow $p$-subgroup $P$ is $p$-nilpotent if $P \le Z(N_G(P))$, or equivalently if $P$ is abelian and $N_G(P)$ is $p$-nilpotent (see \cite[Theorem 5.13]{Isaacs}). This result has been generalized in many directions, and we shall now consider some of these generalizations. 
	
	First let us introduce some notation. Let $P$ be a $p$-group and $i$ be a positive integer. Then the subgroup of $P$ generated by all elements $x$ of $P$ with $x^{p^i} = 1$ is denoted by $\Omega_i(P)$. We set $\Omega(P) := \Omega_1(P)$ if $p$ is odd and $\Omega(P) := \Omega_2(P)$ if $p = 2$. 
	
	In 1974, Laffey \cite{Laffey} published the following generalization of Burnside's $p$-nilpotency criterion: If $G$ is a group and $P$ is a Sylow $p$-subgroup of $G$, then $G$ is $p$-nilpotent if $\Omega(P) \le Z(P)$ and $N_G(P)$ is $p$-nilpotent (see \cite[p. 136]{Laffey}). 
	
	In 2000, Ballester-Bolinches and Guo \cite{BBG} published the following more general result: If $G$ is a group and $P$ is a Sylow $p$-subgroup of $G$, then $G$ is $p$-nilpotent if $\Omega(P \cap G') \le Z(N_G(P))$ (see \cite[Theorem 1]{BBG}). Moreover, they proved that a group $G$ with Sylow $2$-subgroup $P$ is $2$-nilpotent if $\Omega_1(P \cap G') \le Z(P)$, $P$ is quaternion-free and $N_G(P)$ is $2$-nilpotent (see \cite[Theorem 2]{BBG}). Here, a group is said to be quaternion-free if it has no section isomorphic to the quaternion group of order $8$. 
	
	For a non-empty formation $\mathfrak{F}$ and a group $G$, we use $G^{\mathfrak{F}}$ to denote the \textit{$\mathfrak{F}$-residual} of $G$, i.e. $G^{\mathfrak{F}}$ is the smallest normal subgroup of $G$ whose quotient lies in $\mathfrak{F}$. As usual, $\mathfrak{N}$ denotes the formation of all nilpotent groups. 
	
	In 2004, Asaad \cite{Asaad} published the following result, which extends the above-mentioned results of Ballester-Bolinches and Guo.

	\begin{customthm}{1}\label{one}
		(\cite[Theorem 1 (a) $\Leftrightarrow$ (b)]{Asaad}) {\it Let $P$ be a Sylow $p$-subgroup of a group $G$. If $p = 2$, assume that $P$ is quaternion-free. Then the following statements are equivalent:
			\begin{enumerate}
				\item[\textnormal{(1)}] $G$ is $p$-nilpotent. 
				\item[\textnormal{(2)}] $N_G(P)$ is $p$-nilpotent and $\Omega_1(G^{\mathfrak{N}} \cap P \cap P^x) \le Z(P)$ for all $x \in G \setminus N_G(P)$.
			\end{enumerate}}
	\end{customthm} 
	
	Recall that a group $G$ is said to be \textit{$p$-solvable} if any chief factor of $G$ is either a $p$-group or a $p'$-group. A $p$-solvable group $G$ is called \textit{$p$-supersolvable} if every $p$-chief factor of $G$ has order $p$. The formation of all $p$-supersolvable groups is denoted by $\mathfrak{U}_p$. Note that every $p$-nilpotent group is $p$-supersolvable. In particular, we have $G^{\mathfrak{U}_p} \le G^{\mathfrak{N}_p}$ for any group $G$. 
	
	By \cite[Theorem 2.2]{WWL}, if $P$ is a Sylow $p$-subgroup of a group $G$, then we have $P \cap G^{\mathfrak{N}} = P \cap G^{\mathfrak{N}_p}$, where $\mathfrak{N}_p$ denotes the formation of all $p$-nilpotent groups. Therefore, the nilpotent residual $G^{\mathfrak{N}}$ can be replaced by the $p$-nilpotent residual $G^{\mathfrak{N}_p}$ in Theorem \ref{one}. In view of this observation, it is natural to ask whether Theorem \ref{one} still remains true when the nilpotent residual $G^{\mathfrak{N}}$ is replaced by the $p$-supersolvable residual $G^{\mathfrak{U}_p}$. We will show that $G^{\mathfrak{N}}$ cannot only be replaced by $G^{\mathfrak{U}_p}$, but that it is enough to require in (2) that the normalizer $N_G(P)$ is $p$-nilpotent and that every minimal subgroup of $G^{\mathfrak{U}_p} \cap P \cap P^x$ is central in $P$ or complemented in $G$ for all $x \in G \setminus N_G(P)$. We will also show that Theorem \ref{one} remains true if we additionally replace "$p$-nilpotent" by "$p$-supersolvable" in both statements of the theorem and furthermore assume $G$ to be $p$-solvable. 
	
	In fact, our results are slightly more general than just stated, and they also deal with the case that $P \in \mathrm{Syl}_p(G)$ is not quaternion-free when $p = 2$. In order to state our results in full generality, we introduce the following definition. 
	
	\begin{customdef}{2}\label{two}
	Let $G$ be a group and $H \le K \le G$. Then $H$ is said to be \textit{$c$-embedded} in $G$ with respect to $K$ if there is a subgroup $B$ of $G$ such that $G = HB$ and $H \cap B \le Z(K)$.	
	\end{customdef} 
	
	Let $G$ be a group and $H \le K \le G$ such that $H$ is $c$-embedded in $G$ with respect to $K$, so that there exists a subgroup $B$ of $G$ with $G = HB$ and $H \cap B \le Z(K)$. If $H \cap B = 1$ then $H$ is complemented in $G$, and if $H \cap B = H$ then $H$ is central in $K$. If $1 < H \cap B < H$, then $H$ can be described as being between complemented in $G$ and central in $K$. Note that if $H$ is a minimal subgroup of $G$, then $H$ is $c$-embedded in $G$ with respect to $K$ if and only if $H$ is complemented in $G$ or central in $K$. 
	
	Having introduced the concept of $c$-embedded subgroups, we can now state our main results.
	
	\begin{customthm}{A}\label{A}
	{\it Let $P$ be a Sylow $p$-subgroup of a $p$-solvable group $G$. Then $G$ is $p$-supersolvable if and only if the following two conditions are satisfied: 
	\begin{enumerate}
		\item[\textnormal{(1)}] $N_G(P)$ is $p$-supersolvable.  
		\item[\textnormal{(2)}] For all $x \in G \setminus N_G(P)$, the following hold: Any subgroup of $G^{\mathfrak{U}_p} \cap P \cap P^x$ with order $p$ is $c$-embedded in $G$ with respect to $P$. Moreover, if $p = 2$ and $G^{\mathfrak{U}_p} \cap P \cap P^x$ is not quaternion-free, any cyclic subgroup of $G^{\mathfrak{U}_p} \cap P \cap P^x$ with order $4$ is $c$-embedded in $G$ with respect to $P$. 
	\end{enumerate}}	
	\end{customthm}  
	
	The $p$-solvability condition on $G$ cannot be dropped in Theorem \ref{A}. For example, the conditions (1) and (2) from Theorem \ref{A} are satisfied for $G = A_5$ and $p = 5$ since a Sylow $5$-subgroup of $A_5$ is cyclic of order $5$ and has a normalizer of order $10$, but $A_5$ is not $5$-supersolvable. 
	
	\begin{customthm}{B}\label{B}
		{\it Let $P$ be a Sylow $p$-subgroup of a group $G$. Then $G$ is $p$-nilpotent if and only if the following two conditions are satisfied: 
			\begin{enumerate}
				\item[\textnormal{(1)}] $N_G(P)$ is $p$-nilpotent.  
				\item[\textnormal{(2)}] For all $x \in G \setminus N_G(P)$, the following hold: Any subgroup of $G^{\mathfrak{U}_p} \cap P \cap P^x$ with order $p$ is $c$-embedded in $G$ with respect to $P$. Moreover, if $p = 2$ and $G^{\mathfrak{U}_p} \cap P \cap P^x$ is not quaternion-free, any cyclic subgroup of $G^{\mathfrak{U}_p} \cap P \cap P^x$ with order $4$ is $c$-embedded in $G$ with respect to $P$. 
		\end{enumerate}}	
	\end{customthm} 
	
	\section{Preliminaries}
	In this section, we collect some lemmas needed for the proofs of our main results. The following well-known result can be deduced from \cite[Theorem 1 and Proposition 1]{BP}.
	
	\begin{customlem}{3}\label{three}
		{\it Let $G$ be a $p$-solvable minimal non-$p$-supersolvable group. Then the following hold: 
			\begin{enumerate}
				\item[\textnormal{(1)}] $G^{\mathfrak{U}_p}$ is a $p$-group, $G^{\mathfrak{U}_p}$ has exponent $p$ if $p$ is odd and exponent at most $4$ if $p = 2$. 
				\item[\textnormal{(2)}] $G^{\mathfrak{U}_p}/\Phi(G^{\mathfrak{U}_p})$ is a chief factor of $G$. 
		\end{enumerate}}
	\end{customlem} 

\begin{customlem}{4}\label{four}
	{\it Let $G$ be a $p$-solvable minimal non-$p$-supersolvable group, $H$ be a proper subgroup of $G^{\mathfrak{U}_p}$ and $B$ be a subgroup of $G$ such that $G = HB$. Then $B = G$.}
\end{customlem} 

\begin{proof}
	Assume for the sake of contradiction that $B$ is a proper subgroup of $G$. Let $M$ be a maximal subgroup of $G$ such that $B \le M$. Then $G^{\mathfrak{U}_p} \not\le M$ because otherwise $G = HB \le M$. Thus $G = MG^{\mathfrak{U}_p}$.
	
	We have $\Phi(G^{\mathfrak{U}_p}) \le G^{\mathfrak{U}_p} \cap \Phi(G) \le G^{\mathfrak{U}_p} \cap M$. By Lemma \ref{three}(1), $G^{\mathfrak{U}_p}$ is a $p$-group. Since $G^{\mathfrak{U}_p}/\Phi(G^{\mathfrak{U}_p})$ is elementary abelian, $(G^{\mathfrak{U}_p} \cap M)/\Phi(G^{\mathfrak{U}_p})$ is normal in $G^{\mathfrak{U}_p}/\Phi(G^{\mathfrak{U}_p})$. So $G^{\mathfrak{U}_p} \cap M$ is normal in $G^{\mathfrak{U}_p}$. Since $G^{\mathfrak{U}_p} \cap M$ is also normal in $M$ and since $G = MG^{\mathfrak{U}_p}$, it follows that $G^{\mathfrak{U}_p} \cap M \trianglelefteq G$.
	
	As $G^{\mathfrak{U}_p}/\Phi(G^{\mathfrak{U}_p})$ is a chief factor of $G$ by Lemma \ref{three}(2), it follows that either $G^{\mathfrak{U}_p} \cap M = \Phi(G^{\mathfrak{U}_p})$ or $G^{\mathfrak{U}_p} \cap M = G^{\mathfrak{U}_p}$. In the latter case, we have $H \le M$ and so $G = HB \le M$, a contradiction. Thus $G^{\mathfrak{U}_p} \cap M = \Phi(G^{\mathfrak{U}_p})$.
	
	Now we have $G^{\mathfrak{U}_p} = G^{\mathfrak{U}_p} \cap HB = H(G^{\mathfrak{U}_p} \cap B) \le H(G^{\mathfrak{U}_p} \cap M) = H\Phi(G^{\mathfrak{U}_p})$. Thus $G^{\mathfrak{U}_p} = H\Phi(G^{\mathfrak{U}_p})$ and so $G^{\mathfrak{U}_p} = H$. This is a contradiction since $H$ is assumed to be a proper subgroup of $G^{\mathfrak{U}_p}$. So we have $B = G$.
\end{proof}

\begin{customlem}{5}\label{five}
	{\it Let $G$ be a group and $L$ be a subgroup of $G$. Then $L^{\mathfrak{U}_p} \le G^{\mathfrak{U}_p}$. 
	}
\end{customlem}

\begin{proof}
	Set $L_0 := L \cap G^{\mathfrak{U}_p}$. Then $L/L_0 \cong LG^{\mathfrak{U}_p}/G^{\mathfrak{U}_p} \le G/G^{\mathfrak{U}_p}$. As subgroups of $p$-supersolvable groups are $p$-supersolvable, it follows that $L/L_0$ is $p$-supersolvable. Thus $L^{\mathfrak{U}_p} \le L_0 \le G^{\mathfrak{U}_p}$.
\end{proof}

\begin{customlem}{6}\label{six}
{\it Let $G$ be a $p$-supersolvable group and $P$ be a Sylow $p$-subgroup of $G$. Suppose that $N_G(P)$ is $p$-nilpotent. Then $G$ is $p$-nilpotent.  
}
\end{customlem}

\begin{proof}
Set $\overline G := G/O_{p'}(G)$. Then $\overline G$ is $p$-supersolvable and $O_{p'}(\overline G) = 1$. So $\overline G$ is $p$-closed by \cite[Lemma 2.1.6]{BEA}. Consequently $\overline P \trianglelefteq \overline G$, and so $\overline G = N_{\overline G}(\overline P) = \overline{N_G(P)}$. Now, since $N_G(P)$ is $p$-nilpotent, we have that $\overline G = \overline{N_G(P)}$ is $p$-nilpotent. This implies that $G$ is $p$-nilpotent.
\end{proof} 

\section{Proofs of Theorems A and B}
\begin{proof}[Proof of Theorem \ref{A}] Let $G$ and $P$ be as in the statement of Theorem \ref{A}. Assume that $G$ is $p$-supersolvable. Then, since subgroups of $p$-supersolvable groups are $p$-supersolvable, we have that $N_G(P)$ is $p$-supersolvable, whence condition (1) from Theorem \ref{A} is satisfied. Also $G^{\mathfrak{U}_p} = 1$, so that condition (2) from Theorem \ref{A} is trivially satisfied. 
	
	Suppose now that conditions (1) and (2) from Theorem \ref{A} are satisfied. We have to show that $G$ is $p$-supersolvable. To prove this, we assume that $G$ is not $p$-supersolvable, and we are going to derive a contradiction from this assumption. Since $G$ is not $p$-supersolvable, $G$ has a minimal non-$p$-supersolvable subgroup, say $L$. Without loss of generality, we assume that $P \cap L \in \mathrm{Syl}_p(L)$. As $G$ is $p$-solvable, we have that $L$ is $p$-solvable, and so $L^{\mathfrak{U}_p}$ is a $p$-group by Lemma \ref{three}(1). In particular $L^{\mathfrak{U}_p} \le P \cap L$. In order to obtain the desired contradiction, we proceed in a number of steps. 
	
	\medskip
	1) \textit{Any subgroup of $L^{\mathfrak{U}_p}$ with order $p$ is $c$-embedded in $G$ with respect to $P$. Moreover, if $p = 2$ and $L^{\mathfrak{U}_p}$ is not quaternion-free, any cyclic subgroup of $L^{\mathfrak{U}_p}$ with order $4$ is $c$-embedded in $G$ with respect to $P$.}
	
	Since $N_G(P)$ is $p$-supersolvable and $L$ is not $p$-supersolvable, we have $L \not\le N_G(P)$. Let $x \in L \setminus N_G(P)$. Using Lemma \ref{five} and the fact that $L^{\mathfrak{U}_p} \trianglelefteq L$, we see that $L^{\mathfrak{U}_p} = L^{\mathfrak{U}_p} \cap (L^{\mathfrak{U}_p})^x \cap G^{\mathfrak{U}_p} \le P \cap P^x \cap G^{\mathfrak{U}_p}$. Since condition (2) from Theorem \ref{A} is satisfied by assumption, it follows that any subgroup of $L^{\mathfrak{U}_p}$ with order $p$ is $c$-embedded in $G$ with respect to $P$. 
	
	Assume that $p = 2$ and that $L^{\mathfrak{U}_p}$ is not quaternion-free. Then $P \cap P^x \cap G^{\mathfrak{U}_p}$ is not quaternion-free either, and so the validity of condition (2) from Theorem \ref{A} implies that any cyclic subgroup of $L^{\mathfrak{U}_p}$ with order $4$ is $c$-embedded in $G$ with respect to $P$. 
	
	\medskip
	2) \textit{If $H$ is a proper subgroup of $L^{\mathfrak{U}_p}$ which is $c$-embedded in $G$ with respect to $P$, then $H \le Z(P)$.}
	
	Let $H$ be a proper subgroup of $L^{\mathfrak{U}_p}$ such that $H$ is $c$-embedded in $G$ with respect to $P$. Hence there is a subgroup $B$ of $G$ such that $G = HB$ and $H \cap B \le Z(P)$. Set $B_0 := L \cap B$. Then $L = L \cap HB = HB_0$. Lemma \ref{four} implies that $L = B_0 \le B$. So it follows that $H = H \cap B \le Z(P)$. 
	
	\medskip
	3) \textit{$L^{\mathfrak{U}_p} \le Z(P)$.} 
	
	Suppose that $L^{\mathfrak{U}_p}$ has exponent $p$. Let $x \in L^{\mathfrak{U}_p}$. We show that $x \in Z(P)$. Clearly, we only need to consider the case $x \ne 1$. Then $|\langle x \rangle| = p$. We have $|L^{\mathfrak{U}_p}| > p$ since $L$ would be $p$-supersolvable otherwise. Hence $\langle x \rangle$ is a proper subgroup of $L^{\mathfrak{U}_p}$. So 1) and 2) imply that $x \in Z(P)$. As $x$ was arbitrarily chosen, it follows that $L^{\mathfrak{U}_p} \le Z(P)$. 
	
	Suppose now that $L^{\mathfrak{U}_p}$ does not have exponent $p$. Then, by Lemma \ref{three}(1), $p = 2$ and $L^{\mathfrak{U}_2}$ has exponent $4$. Therefore $L$ is a group appearing in \cite[Theorem 9]{BE} as a group of Type 3. In particular, $L^{\mathfrak{U}_2} = P \cap L$ is a non-abelian special $2$-group, $\Phi(L^{\mathfrak{U}_2}) \le Z(L)$ and $|L^{\mathfrak{U}_2}/\Phi(L^{\mathfrak{U}_2})| = 2^{2m}$, $|\Phi(L^{\mathfrak{U}_2})| \le 2^m$ for some positive integer $m$.
	
	We claim that $L^{\mathfrak{U}_2}$ is not quaternion-free. Assume that $m = 1$. Then $|L^{\mathfrak{U}_2}| = 8$, and $L^{\mathfrak{U}_2}$ cannot be dihedral because then $L^{\mathfrak{U}_2}/\Phi(L^{\mathfrak{U}_2})$ would not be a chief factor of $L$. Consequently $L^{\mathfrak{U}_2}$ is isomorphic to the quaternion group of order $8$ and in particular not quaternion-free. Assume now that $m > 1$. Let $R$ be a maximal subgroup of $\Phi(L^{\mathfrak{U}_2})$. Then $L^{\mathfrak{U}_2}/R$ is extraspecial of order $2^{2m+1} \ge 2^5$, and so $L^{\mathfrak{U}_2}/R$ has a section isomorphic to the quaternion group of order $8$ (see \cite[Chapter 5, Theorem 5.2]{GO}). Hence $L^{\mathfrak{U}_2}$ is not quaternion-free. 
	
	
	Now let $x \in L^{\mathfrak{U}_2}$. We show that $x \in Z(P)$. Clearly, we only need to consider the case $x \ne 1$. Then $|\langle x \rangle| = 2$ or $4$. Also $\langle x \rangle$ is a proper subgroup of $L^{\mathfrak{U}_2}$. So 1) and 2) imply that $x \in Z(P)$. As $x$ was arbitrarily chosen, it follows that $L^{\mathfrak{U}_2} \le Z(P)$.
	
	\medskip
	4) \textit{The final contradiction.}
	
	Set $N := N_G(L^{\mathfrak{U}_p})$. By 3), we have $P \le C_G(L^{\mathfrak{U}_p})$. Since $C_G(L^{\mathfrak{U}_p}) \trianglelefteq N$, the Frattini argument implies that $N = N_N(P)C_G(L^{\mathfrak{U}_p})$.
	
	Clearly $L^{\mathfrak{U}_p} \trianglelefteq N_N(P)$. Let $1 = L_0 \le L_1 \le \cdots \le L_t = L^{\mathfrak{U}_p}$ be a part of a chief series of $N_N(P)$ below $L^{\mathfrak{U}_p}$, i.e. $L_i \trianglelefteq N_N(P)$ for all $0 \le i \le t$ and $L_{i+1}/L_i$ is a chief factor of $N_N(P)$ for all $0 \le i < t$. Since $N_G(P)$ is $p$-supersolvable and $N_N(P) \le N_G(P)$, we have that $N_N(P)$ is $p$-supersolvable. Consequently $|L_{i+1}/L_i| = p$ for all $0 \le i < t$.
	
	For each $0 \le i \le t$, we have $C_G(L^{\mathfrak{U}_p}) \le C_G(L_i) \le N_G(L_i)$ and hence $N = N_N(P)C_G(L^{\mathfrak{U}_p}) \le N_G(L_i)$. In particular, we have $L_i \trianglelefteq L$ for all $0 \le i \le t$. Consequently, each chief factor of $L$ below $L^{\mathfrak{U}_p}$ has order $p$. So it follows that $L$ is $p$-supersolvable. This contradiction completes the proof. 
\end{proof}

\begin{proof}[Proof of Theorem \ref{B}]
Let $P$ be a Sylow $p$-subgroup of a group $G$. Assume that $G$ is $p$-nilpotent. Then, since subgroups of $p$-nilpotent groups are $p$-nilpotent, we have that $N_G(P)$ is $p$-nilpotent, whence condition (1) from Theorem \ref{B} is satisfied. Also $G^{\mathfrak{U}_p} \le G^{\mathfrak{N}_p} = 1$, so that condition (2) from Theorem \ref{B} is trivially satisfied. 

Let $\mathfrak{Y}$ denote the class of all groups $G$ such that condition (2) from Theorem \ref{B} is satisfied for any Sylow $p$-subgroup $P$ of $G$. Note that if condition (2) from Theorem \ref{B} is satisfied for one Sylow $p$-subgroup of a group $G$, then it is satisfied for all Sylow $p$-subgroups of $G$, so that $G \in \mathfrak{Y}$. Let $\mathscr{Z}_p$ denote the class of all $\mathfrak{N}$-groups with $p$-nilpotent normalizers of Sylow $p$-subgroups. To complete the proof of Theorem \ref{B}, we show that the class $\mathscr{Z}_p$ is contained in the class $\mathfrak{N}_p$ of all $p$-nilpotent groups. Suppose that this is not true, and choose a non-$p$-nilpotent $\mathscr{Z}_p$-group $G$ of minimal order.

Arguing similarly as in \cite[Example 2]{BGLS}, we see that $\mathfrak{Y}$ is subgroup-closed and that $X/N \in \mathfrak{Y}$ whenever $X \in \mathfrak{Y}$ and $N$ is a normal $p'$-subgroup of $X$. Applying \cite[Theorem A]{BGLS}, we conclude that $G$ is $p$-solvable. 

Let $P \in \mathrm{Syl}_p(G)$. Then $N_G(P)$ is $p$-nilpotent and hence $p$-supersolvable as $G \in \mathscr{Z}_p$, whence condition (1) from Theorem \ref{A} is satisfied. Also, condition (2) from Theorem \ref{A} is satisfied since it is identical to condition (2) from Theorem \ref{B}, which holds as $G \in \mathfrak{Y}$. So Theorem \ref{A} implies that $G$ is $p$-supersolvable. Applying Lemma \ref{six}, we conclude that $G$ is $p$-nilpotent. This contradiction shows that $\mathscr{Z}_p$ is contained in $\mathfrak{N}_p$, as wanted.	
\end{proof}
	
\section{Remarks and open questions}
One might wonder whether Theorems \ref{A} and \ref{B} remain true when, in condition (2) of the theorems, ``$c$-embedded in $G$ with respect to $P$'' is replaced by ``$c$-embedded in $P$ with respect to $P$''. For the case $p = 2$, the answer is negative, as the following example shows. 

\begin{customex}{7}
Let $G := S_4$, the symmetric group of degree $4$, and let $P$ be a Sylow $2$-subgroup of $G$. Then $P$ is dihedral of order $8$, and we have $N_G(P) = P$. Any subgroup of $P$ with order $2$ is either complemented or central in $P$ and thus $c$-embedded in $P$ with respect to $P$. However, $G$ is not $2$-nilpotent (or, equivalently, not $2$-supersolvable).	
\end{customex} 

For the case $p = 3$, the following example shows that Theorem \ref{A} does not remain true when ``$c$-embedded in $G$ with respect to $P$'' is replaced by ``$c$-embedded in $P$ with respect to $P$''. 

\begin{customex}{8}
	Let $G$ be the group indexed in GAP \cite{GAP} as SmallGroup(216,153). Then $G$ is solvable and hence $3$-solvable. Let $P$ be a Sylow $3$-subgroup of $G$. Then $N_G(P)$ is $3$-supersolvable, and any subgroup of $P$ with order $3$ is complemented or central in $P$ and thus $c$-embedded in $P$ with respect to $P$. However, $G$ is not $3$-supersolvable.
\end{customex} 

We were not able to answer the following question.

\begin{customqu}{9}\label{nine}
Suppose that $p$ is odd. Does Theorem \ref{B} remain true when, in condition (2) from Theorem B, ``$c$-embedded in $G$ with respect to $P$'' is replaced by ``$c$-embedded in $P$ with respect to $P$''?
\end{customqu}

Using GAP \cite{GAP}, we have checked that there are no counterexamples of order up to $2000$. 

Wei, Wang and Liu \cite{WWL} obtained the following characterization of $p$-nilpotent groups: A group $G$ with Sylow $p$-subgroup $P$ is $p$-nilpotent if and only if every minimal subgroup of $P \cap G^{\mathfrak{N}_p}$ is complemented in $P$ and $N_G(P)$ is $p$-nilpotent (see \cite[Corollary 2.3]{WWL}). In a less general form, this had been proved before by Guo and Shum \cite[Theorem 2.1]{GS}. Note that if the correct answer to Question \ref{nine} is positive, then this would generalize both Theorem \ref{B} and the mentioned result of Wei, Wang and Liu. 

\section*{Acknowledgements}
I am grateful to Professor Adolfo Ballester-Bolinches for drawing my attention to the papers \cite{BE} and \cite{BGLS}.

\end{document}